\documentclass[conference]{IEEEtran}
\usepackage{cite}
\usepackage{amsmath,amssymb,amsfonts,amsthm}
\usepackage{enumitem}
\usepackage{mathtools}
\usepackage{amsbsy}
\usepackage{mathrsfs}
\usepackage[active]{srcltx}
\usepackage{bm}
\usepackage{fixmath}
\usepackage{enumitem}
\usepackage{textcomp}
\def\BibTeX{{\rm B\kern-.05em{\sc i\kern-.025em b}\kern-.08em
    T\kern-.1667em\lower.7ex\hbox{E}\kern-.125emX}}

\newtheorem{theorem}{Theorem}

\newtheorem{remark}{Remark}

\begin{document}

\title{A bipolar Hardy inequality on Finsler manifolds}

\author{\IEEEauthorblockN{Mester \'{A}gnes}
\IEEEauthorblockA{\textit{Doctoral School of Applied Mathematics} \\
\textit{\'Obuda University}\\
Budapest, Hungary \\
mester.agnes@yahoo.com}
\and
\IEEEauthorblockN{Krist\'aly Alexandru}
\IEEEauthorblockA{\textit{Institute of Applied Mathematics} \\
\textit{\'Obuda University}\\
Budapest, Hungary \\
and\\
\textit{Department of Economics}\\
\textit{ Babe\c s-Bolyai University}\\
Cluj-Napoca, Romania\\
kristaly.alexandru@nik.uni-obuda.hu;\\
alex.kristaly@econ.ubbcluj.ro}
}

\maketitle

\begin{abstract}
We establish a bipolar Hardy inequality on comp\-lete, not necessarily reversible Finsler manifolds. 
We show that our result strongly depends on the geometry of the Finsler structure, namely on the reversibility constant $r_F$ and the uniformity constant $l_F$.
Our result represents a Finslerian counterpart of the Euclidean multipolar Hardy inequality due to Cazacu and Zuazua \cite{CZ} and the Riemannian case considered by  Faraci, Farkas and Krist\'aly \cite{FFK}.
\end{abstract}


\section{Introduction and main results}

The classical Hardy inequality states that
\begin{equation} \label{unipolar_Hardy}
 \int_{\mathbb{R}^n} |\nabla u|^2 dx  \geq  \frac{(n-2)^2}{4}  \int_{\mathbb{R}^n} \frac{u^2}{|x|^2} dx  , 
\quad \forall u \in C_0^{\infty}(\mathbb{R}^n) ,
\end{equation}
where the constant $\frac{(n-2)^2}{4}$ is optimal and not achieved, see Hardy, Littlewood and P\'olya \cite{HLP}.

A challenging direction of extension consists of the study of multipolar Hardy inequalities, 
motivated by the applications in molecular physics, quantum cosmology and combustion models, 
see Bosi, Dolbeault and Esteban \cite{BDE}, Felli, Marchini and Terracini \cite{FMT}, Guo, Han and Niu \cite{GHN} and references therein. 

The optimal multipolar extension of the unipolar inequa\-lity \eqref{unipolar_Hardy} was proved by Cazacu and Zuazua \cite{CZ}: 
\begin{align} \label{multipolar_Hardy_flat}
& \int_{\mathbb{R}^n} |\nabla u|^2 dx  \geq  \nonumber \\
& \geq \frac{(n-2)^2}{m^2} \sum_{1\leq i<j\leq m} \int_{\mathbb{R}^n} \bigg| \frac{x-x_i}{|x-x_i|^2} - \frac{x-x_j}{|x-x_j|^2} \bigg|^2 u^2 dx  , 
\end{align}
$\forall u \in C_0^{\infty}(\mathbb{R}^n)$, 
where $x_1, \dots, x_m \in \mathbb{R}^n$ represent pairwise distinct poles, $m \geq 2$, $n \geq 3$, and the constant $\frac{(n-2)^2}{m^2}$ is sharp.

Recently, there has been a growing attempt to develop the theory of Hardy inequalities on Riemannian and Finsler manifolds, see e.g. Kombe and \"{O}zaydin \cite{KO}, D'Ambrosio and Dipierro \cite{DD}, Xia \cite{X}, Yang, Su and Kong \cite{YSK}, Farkas, Krist\'aly and Varga \cite{FKV}, Krist\'aly and Repov\v{s} \cite{KR}, and Yuan, Zhao and Shen \cite{YZS}.

In 2018, Faraci, Farkas and Krist\'aly \cite{FFK} proved multipolar Hardy inequalities on complete Riemannian manifolds, obtaining the curved analogue of inequality \eqref{multipolar_Hardy_flat}. In order to present this result, let us consider an $n$-dimensional complete Riemannian manifold $(M,g)$ with $n \geq 3$. 
Let $dv_g$ and $d_g: M \times M \to [0,\infty)$ denote the canonical volume form and the distance function defined on $M$, 
induced by the Riemannian metric $g$.
Furthermore, $\nabla_g$ and $\Delta_g$ stand for the gradient operator and Laplace-Beltrami operator defined on $(M,g)$. 
For the sake of brevity, in the sequel let $|\cdot|$ denote the norm associated with the Riemannian metric $g$.  
Finally, let $x_1, \dots, x_m \in M$ be the set of pairwise distinct poles, where $m \geq 2$. Then the following multipolar Hardy inequality holds:
\begin{align} \label{multipolar_Hardy_curved}
& \int_{M} |\nabla_g u|^2 dv_g  \geq  \nonumber \\
& \geq \frac{(n-2)^2}{m^2} \sum_{1\leq i<j\leq m} \int_{M} \bigg| \frac{\nabla_g d_i}{d_i} - \frac{\nabla_g d_j}{d_j} \bigg|^2 u^2 dv_g \nonumber \\
& + \frac{n-2}{m} \sum_{i=1}^m \int_{M} \frac{d_i \Delta_g d_i - (n-1)}{d_i^2}  u^2 dv_g  , 
\end{align}
$\forall u \in C_0^{\infty}(M)$, 
where $d_i = d_g( x_i, \cdot)$ denotes the Riemannian distance from the pole $x_i \in M$, $i = \overline{1,m}$. 
The constant $\frac{(n-2)^2}{m^2}$ is sharp in the bipolar case, i.e. when $m=2$. 

A few remarks are in order considering the importance of the last term in inequality \eqref{multipolar_Hardy_curved} (for the full discussion see \cite{FFK}):
\begin{itemize} 
	\item
	if the Ricci curvature of the manifold satisfies $Ric(M,g) \geq c_0 (n-1)g$ for some $c_0 > 0$, then it can be proven that the last term is negative, thus modifying the analogue of the flat case \eqref{multipolar_Hardy_flat} in order to hold true.
	\item
	in the negatively curved case, by using a suitable Laplace comparison theorem (see Wu and Xin \cite{WX}), one can  prove that the last term in \eqref{multipolar_Hardy_curved} enables us to obtain stronger inequality when stronger curvature is assumed. 
	\item
	if $(M,g) = (\mathbb{R}^n, g_0)$ is the standard Euclidean space, then $d_i(x) = |x - x_i|, \forall x \in \mathbb{R}^n$, $|\cdot|$ being the Euclidean norm, thus the last term vanishes, and we obtain  \eqref{multipolar_Hardy_flat}.     
\end{itemize}

The purpose of this paper is to study multipolar Hardy inequalities on complete, not necessarily reversible Finsler manifolds. We notice that the obtained results heavily depend on the non-Riemannian nature of Finsler structures, expressed in terms of the reversibility constant $r_F$ and uniformity constant $l_F$.

In order to present our results, let $(M,F)$ be a complete Finsler manifold, and let us denote by $\mathrm{div}$, $\nabla_F$ and $\Delta_F$ the divergence, gradient and Finsler-Laplace operator determined by the Finsler structure $F$. 

Furthermore, $dv_F$ and $d_F: M \times M \to [0,\infty)$ denote the Busemann-Hausdorff volume form and distance function defined on $(M, F)$, respectively, while $F^{*}$ is the polar transform of $F$ and $J^*:T^*M \to TM$ is the Legendre transform. Finally, let $r_F \in [1, \infty)$ and $l_F \in (0,1]$ be the reversibility and uniformity constant of the Finsler manifold $(M,F)$ (for the detailed definitions see Section \ref{preliminaries}). 
Our first result reads as follows:

\begin{theorem} \label{thm_multipolar_Hardy_Finsler}

Let $(M,F)$ be a complete $n$-dimensional Finsler manifold with $n \geq 3$ and $l_F > 0$, 
and consider the set of pairwise distinct poles $\{x_1, \dots, x_m\} \subset M$, where $m \geq 2$. 
Then 
\begin{align} \label{multipolar_Hardy_Finsler}
& \left( 2 - \frac{l_F^2}{r_F^2}\right) \int_M F^{*2}( Du ) dv_F   \geq \nonumber \\
& \geq  (l_F-2)\frac{(n-2)^2}{m^2} \int_M F^{*2} \Big( \sum_{i=1}^m \frac{Dd_i}{d_i} \Big) u^2 dv_F   \nonumber \\
& + l_F \frac{n-2}{m} \int_M \mathrm{div} \Big(J^* \Big(\sum_{i=1}^m \frac{Dd_i}{d_i}\Big) \Big) u^2  dv_F 
\end{align}
holds for every nonnegative function $u \in C_0^{\infty}(M)$,
where $d_i(x) = d_F(x, x_i)$ denotes the Finslerian distance from the point $x$ to the pole $x_i$, $i = \overline{1,m}$. 
\end{theorem}

We shall prove in Section \ref{proof} that when $(M,F) = (M,g)$ is a Riemannian manifold, 
then the inequality above is equivalent with \eqref{multipolar_Hardy_curved}, 
meaning that our result extends the multipolar Hardy inequality 
obtained by \cite[Theorem 1.1]{FFK} to the case of complete Finsler manifolds.

By using Theorem \ref{thm_multipolar_Hardy_Finsler} in the case $m=2$, we obtain a bipolar Hardy inequality.

\begin{theorem} \label{thm_bipolar_Hardy_Finsler}
Let $(M,F)$ be a complete $n$-dimensional Finsler manifold with $n \geq 3$ and $l_F > 0$. 
Let $x_1, x_2 \in M, x_1 \neq x_2$ be two poles. 
Then 
\begin{align}	\label{bipolar_Hardy_Finsler}
&\int_M F^{*2}( Du ) dv_F   \geq \nonumber \\ 
& \geq \frac{l_F(2-l_F)}{ 2 - \big(\frac{l_F}{r_F}\big)^2 } \frac{(n-2)^2}{4} 
\int_M  F^{*2} \Big(\frac{Dd_2}{d_2} - \frac{Dd_1}{d_1} \Big) u^2 dv_F  \nonumber \\
& + \frac{l_F}{ 2 - \big(\frac{l_F}{r_F}\big)^2 } \frac{n-2}{2} \int_M  \mathrm{div} \bigg(J^*\Big(\frac{Dd_1}{d_1} + \frac{Dd_2}{d_2} \Big) \bigg) u^2 dv_F  \nonumber   \\
& -  \frac{2-l_F}{ 2 - \big(\frac{l_F}{r_F}\big)^2 } \frac{(n-2)^2}{2} \int_M \Big( \frac{1}{d_1^2} + \frac{1}{d_2^2} \Big) u^2 dv_F
\end{align}
holds for every nonnegative function $u \in C_0^{\infty}(M)$. 
\end{theorem}

These results seem to be the first contributions considering multipolar Hardy inequalities in the Finslerian setting.

The next section recalls the notions of Finsler geometry necessary for our further developments. 
Section \ref{proof} contains the proofs of Theorems \ref{thm_multipolar_Hardy_Finsler} and \ref{thm_bipolar_Hardy_Finsler}, 
as well as the proof of equivalence between inequalities \eqref{multipolar_Hardy_curved} and \eqref{multipolar_Hardy_Finsler} in the Riemannian setting.

\section{Elements of Finsler geometry}    \label{preliminaries}

In this section we recall several notions from Finsler geo\-metry, see Bao, Chern and Shen \cite{BCS},  Farkas,  Krist\'aly and  Varga \cite{FKV} and Ohta and Sturm \cite{OS}.

\subsection{Finsler structure, Chern connection, completeness} 

Let $M$ be a connected $n$-dimensional differentiable manifold, and $TM=\bigcup_{x \in M}T_{x} M $ its tangent bundle, 
where $T_{x} M$ denotes the tangent space at the point $x \in M$.

The pair $(M,F)$ is called a Finsler manifold, if $F: TM \to [0,\infty)$ is a continuous function satisfying the following conditions:
\begin{enumerate}[label=(\roman*)]
\item $F$ is of class $C^{\infty}$ on the set $TM \setminus \{ 0 \}$; 
\item $F(x,\lambda y) = \lambda F(x,y)$, for every $\lambda \geq 0$ and $(x,y) \in TM$;
\item the Hessian matrix $\Big( g_{ij}(x,y) \Big) = \left( \left[ \frac{1}{2}F^{2}(x,y)\right] _{y^{i}y^{j}} \right)$
is positive definite for every $(x,y)\in TM \setminus \{0\}.$
\end{enumerate}

The function $F$ is called the Finsler structure on $M$. 
If, in addition, $F(x,\lambda y) = |\lambda| F(x,y)$ holds for all $\lambda \in \mathbb{R}$ and $(x,y) \in TM$, then the Finsler manifold is called reversible. Otherwise, $(M,F)$ is said to be nonreversible.

Let $T^*M = \bigcup_{x \in M}T^*_{x} M $ denote the cotangent bundle of $M$, where $T^*_{x} M$ is the dual space of $T_{x} M$. In the following, for every point $x \in M$ let $( \frac{\partial}{\partial x^i} )_{i=\overline{1,n}}$ be the canonical basis of the tangent space $T_xM$, and $( d x^i )_{i=\overline{1,n}}$ be the dual basis of $T^*_xM$, where $(x^i)_{i=\overline{1,n}}$ is a local coordinate system.

Now let $\pi^* TM$ be the pull-back tangent bundle of $TM$, induced by the natural projection $\pi: TM \setminus \{0\} \to M$, see Bao, Chern and Shen \cite[Chapter~2]{BCS}. 
Thus $\pi^* TM$ is the collection of all pairs $(v;w)$ with $v = (x,y) \in TM \setminus \{0\}$ and $w \in T_xM$. The pull-back tangent bundle admits a natural local basis defined by $\partial_i|_v = (v; \frac{\partial}{\partial x^i})$, and a natural Riemannian metric induced by the Hessian matrices $(g_{ij})$, i.e.
$$ g_{(x,y)}(\partial_i|_v, \partial_j|_v) =  g_{ij}(x,y) .$$
The metric $g$ is called the fundamental tensor on $\pi^* TM$.  

Unlike the Riemannian metric, the Finsler structure $F$ does not induce a unique natural connection on the Finsler manifold $(M, F)$. 
However, on the pull-back tangent bundle $\pi^* TM$ it is possible to define a linear, torsion-free and almost metric-compatible connection called the Chern connection, see Bao, Chern and Shen \cite[Chapter~2]{BCS}.   
The Chern connection induces the notions of covariant derivative and  parallelism of a vector field along a curve. For example, let us denote by $D_yV$ the covariant derivative of a vector field $V$ in the direction $y \in T_xM$. Then, a vector field $V = V(t)$ is parallel along a curve $\gamma = \gamma(t)$ if $D_{\dot \gamma} V = 0$. 

A curve $\gamma: [a,b] \to M$ is called a geodesic if its velocity field $\dot \gamma$ is parallel along the curve, i.e. if $D_{\dot \gamma} \dot \gamma = 0$.
A Finsler manifold is said to be complete if every geodesic segment $\gamma: [a,b] \to M$ can be extended to a geodesic defined on $\mathbb{R}$.

\subsection{Polar transform and Legendre transform}

Let us consider the polar transform $F^*: T^*M \to [0, \infty)$, which is defined as the dual metric of $F$ on $M$, namely
\begin{equation*}  
F^*(x,\alpha) = \sup_{y \in T_xM \setminus \{0\}} ~ \frac{\alpha(y)}{F(x,y)}.
\end{equation*}

We have that for every $x \in M$, $F^{*2}(x, \cdot)$ is twice differentiable on $T^*_xM \setminus \{0\}$. Thus we can define the dual matrix
$$\Big( g^*_{ij}(x,\alpha) \Big) = \left( \left[ \frac{1}{2}F^{*2}(x,\alpha)\right] _{\alpha^{i}\alpha^{j}} \right) ,$$
for every  $\alpha = \sum_{i=1}^n \alpha^i dx^i \in T_x^*M \setminus \{0\}$.

The Legendre transform $J^*:T^*M \to TM$ is defined in the following way: for every $x \in M$ fixed, $J^*$ associates to each $\alpha \in T_x^*M$ the unique maximizer $y \in T_xM$ of the mapping 
$$ y ~ ~ \mapsto ~ ~ \alpha(y) - \frac{1}{2}F^2(x,y) .$$ 
It can be proven that when $J^*(x,\alpha) = (x,y)$, then 
\begin{equation*}  \label{Legendre_transform}
F(x,y) = F^*(x,\alpha) \quad \text{ and } \quad \alpha(y) = F^*(x,\alpha) F(x,y).
\end{equation*}
Moreover, we have the following local characterization of $J^*$. For every $\alpha = \sum_{i=1}^n \alpha^i dx^i \in T_x^*M$, we have that
\begin{equation*}  
J^*(x,\alpha) = \sum_{i=1}^n \frac{\partial}{\partial \alpha_i}\left(\frac{1}{2} F^{*2}(x,\alpha)\right)\frac{\partial}{\partial x^i}  .
\end{equation*}
For further details on the Legendre transform see Bao, Chern and Shen \cite[Section~14.8]{BCS} and Ohta and Sturm \cite{OS}.

\subsection{Hausdorff volume form and distance function}

Let
$ B_x(1) = \Big\{ (y^i) \in \mathbb{R}^n :~ F \Big(x, \sum_{i=1}^n y^i \frac{\partial}{\partial x^i} \Big) < 1 \Big\} \subset \mathbb{R}^n ,$
and define the ratio
$ \sigma_F(x) = \frac{ \omega_n}{\mathrm{Vol}(B_x(1))} ,$
where $\omega_n$ and $\mathrm{Vol}(B_x(1))$ denote the Euclidean volume of the $n$-dimensional unit ball and the set $B_x(1)$, respectively.
The Busemann-Hausdorff volume form is defined as  
\begin{equation*} \label{Hausdorff_measure}
dv_F(x) = \sigma_F(x) d x^1 \land \dots \land d x^n ,
\end{equation*}
see Shen \cite[Section~2.2]{S}. Note that in the following we may omit the parameter $x$ for the sake of brevity.

The distance function $d_F: M \times M \to [0, \infty)$ is defined by
$$ d_F(x_1, x_2) = \inf_\gamma  \int_a^b F( \gamma(t), \dot \gamma(t)) ~ dt,$$
where $\gamma: [a,b] \to M$ is any piecewise differentiable curve 
such that $\gamma(a) = x_1$ and $\gamma(b) = x_2$.
It is immediate that $d_F(x_1, x_2) = 0$ if and only if $x_1 = x_2$ and that $d_F$ verifies the triangle inequality. 
However, $d_F$ is symmetric if and only if $(M,F)$ is a reversible Finsler manifold.

We also recall the eikonal equation, see Shen \cite[Lemma~3.2.3]{S}.
For every point $x_0 \in M$, one has 
\begin{equation} \label{eq_eikonal}
F(x, \nabla_F d_F(x_0, x)) = F^*(x, D d_F(x_0, x)) = 1 
\text{ a.e. } x \in M .
\end{equation}

\subsection{Reversibility and uniformity constants}
   
The reversibility constant of the Finsler manifold $(M,F)$ is defined by 
\begin{equation*} 
r_{F} = \sup_{x\in M} ~ \sup_{\substack{ y\in T_{x} M \setminus \{0\}}} \frac{F(x,y)}{F(x,-y)} ~ \in [1, \infty],
\end{equation*}
measuring how far the Finsler structure $F$ is from being reversible (see Rademacher \cite{R}).
Note that $r_F = 1$ if and only if $(M,F)$ is reversible Finsler manifold.

The uniformity constant of $(M,F)$ is defined by
\begin{equation*}
 l_{F} = \inf_{x\in M} ~ \inf_{y,v,w\in  T_xM\setminus \{0\}}\frac{g_{(x,v)}(y,y)}{g_{(x,w)}(y,y)} ~ \in ~ [0, 1],
\end{equation*}
which measures how much $F$ deviates from being a Riemannian structure.
Indeed, $l_F = 1$ if and only if $(M,F)$ is a Riemannian manifold, see Ohta \cite{O}.

Furthermore, by using the definition of $l_F$, it can be proven that
\begin{align} \label{eq_l_F} 
F^{*2}(x, t \alpha + (1-t)\beta) 
& \leq t F^{*2}(x, \alpha) + (1-t) F^{*2}(x, \beta)  \nonumber \\
& - l_F t(1-t) F^{*2}(x, \beta - \alpha), 
\end{align}
for every $x \in M$, $\alpha, \beta \in T^*_xM$ and $t \in [0,1]$, see Ohta and Sturm \cite{OS}. 

We also have the following implication: if $l_F > 0$ then $r_F < \infty$, see Farkas,  Krist\'aly and  Varga \cite{FKV}.

\subsection{Gradient, divergence, Finsler-Laplace operator}

Let $u: M \to \mathbb{R}$ be a weakly differentiable function. 
Then $Du(x) \in T_x^*M$ denotes the differential of $u$ at every regular point $x \in M$, while the gradient of $u$ at $x$ is defined by
\begin{equation*}  
\nabla_F u(x) = J^*(x, Du(x)).
\end{equation*}
Using the properties of the Legendre transform, it follows that 
\begin{equation*}  \label{metric-co-metric2}
F^*(x, Du(x)) = F(x, \nabla_F u(x)).
\end{equation*}
Also, in local coordinates we can write 
\begin{equation*} 
Du(x) = \sum_{i=1}^n \frac{\partial u}{\partial x^i}(x) dx^i 
\end{equation*}
and
\begin{equation*} 
\nabla_F u(x)=\sum_{i,j=1}^n g_{ij}^*(x,Du(x))\frac{\partial u}{\partial x^i}(x)\frac{\partial}{\partial x^j}.
\end{equation*}
Therefore, the gradient operator $\nabla_F$ is usually nonlinear. 

The divergence operator is defined in a distributional sense, i.e. for every weakly differentiable vector field $V$ on $M$, one has $\mathrm{div} V: M \to \mathbb{R}$ such that 
\begin{equation}  \label{divergence_theorem}
\int_{M} u \mathrm{div} V dv_F = 
-\int_{M} Du (V) dv_F,
\end{equation}
for every $u \in C^{\infty}_0(M)$, see Ohta and Sturm \cite{OS}. 

The Finsler-Laplace operator  
$\Delta_F u = \mathrm{div}(\nabla_F u)$
is defined in a distributional sense as well. 
Note that in general, the Finsler-Laplace operator $\Delta_F$ is nonlinear.

\section{Proof of main results} \label{proof}

In the following let $(M, F)$ be a complete $n$-dimensional Finsler manifold ($n  \geq 3$), such that $l_F > 0$, thus $r_F < \infty$.
We start by proving Theorem \ref{thm_multipolar_Hardy_Finsler}.

\subsection*{\bf Proof of Theorem \ref{thm_multipolar_Hardy_Finsler}.}

For every $x \in M$ and every $\alpha, \beta \in T^*_xM$, we have the following relations:
first, by using \eqref{eq_l_F} for $t=1/2$, one has
\begin{equation} \label{proof1} 
F^{*2}(x, \alpha + \beta) \leq 
2 F^{*2}(x, \alpha) + 2 F^{*2}(x, \beta) - l_F F^{*2}(x, \beta - \alpha).
\end{equation}

Then, due to the strict convexity of $F^{*2}$, we can derive the following inequality:
\begin{equation}  \label{proof2} 
F^{*2}(x, \beta - \alpha) \geq F^{*2}(x, \beta) - 2 \alpha(J^*(x, \beta)) + l_F F^{*2}(x, - \alpha) . 
\end{equation}

Finally, since $r_F < \infty$, we have 
\begin{equation}  \label{proof3} 
F^{*}(x, - \alpha) \geq F^{*}(x, \alpha)/r_F.
\end{equation} 

Using relations \eqref{proof1} -- \eqref{proof3} yields
\begin{align} \label{proof4} 
F^{*2}(x, \alpha + \beta) 
& \leq \left( 2 - \frac{l_F^2}{r_F^2} \right) F^{*2}(x, \alpha) \nonumber \\
& + (2 - l_F) F^{*2}(x, \beta) + 2 l_F ~ \alpha(J^*(x, \beta)) .
\end{align}

Now consider the pairwise distinct poles $x_1, \dots, x_m \in M$ where $m \geq 2$, 
and let $d_i = d_F(\cdot, x_i)$ be the Finslerian distance to the pole $x_i$, $i = \overline{1,m}$. 
Also, let $u \in C_0^{\infty}(M)$ be a function such that $u \geq 0$ on $M$. 
Applying \eqref{proof4} with the choices 
$$ \alpha = Du   ~ \text{ and } ~ \beta = \frac{n-2}{m} u \sum_{i=1}^m \frac{Dd_i}{d_i},$$
then integrating over $M$ results in
\begin{align*} 
0 & \leq \int_M F^{*2} \Big( Du + \frac{n-2}{m}u \sum_{i=1}^m \frac{Dd_i}{d_i} \Big) dv_F  \\ 
  & \leq \left( 2 - \frac{l_F^2}{r_F^2}\right) \int_M F^{*2}( Du ) dv_F    \\
  & +  (2-l_F)\frac{(n-2)^2}{m^2} \int_M F^{*2} \Big( \sum_{i=1}^m \frac{Dd_i}{d_i} \Big) u^2 dv_F \\  
  & +  l_F \frac{n-2}{m} \int_M D(u^2) \Big(J^* \Big(\sum_{i=1}^m \frac{Dd_i}{d_i}\Big) \Big)  dv_F,
\end{align*}
where we omitted the parameter $x$ for the sake of brevity. 
Using the divergence theorem \eqref{divergence_theorem} completes the proof of Theorem \ref{thm_multipolar_Hardy_Finsler}.  \qed

\begin{remark}\rm
Note that if we consider Theorem \ref{thm_multipolar_Hardy_Finsler} in the Riemannian setting, 
then \eqref{multipolar_Hardy_Finsler} becomes equivalent with \eqref{multipolar_Hardy_curved}. 

Indeed, if $(M,F) = (M,g)$ is a Riemannian manifold, then $r_F = l_F = 1$, 
while the operators $\nabla_F$ and $\Delta_F$ coincide with $\nabla_g$ and $\Delta_g$, respectively. 
Moreover, due to the Riesz representation theorem, the tangent space $T_xM$ and its dual space $T_x^*M$ can be identified, 
and the Finsler metrics $F$ and $F^*$ reduce to the norm $|\cdot|$ associated to the Riemannian metric $g$.
Thus the Hardy inequality \eqref{multipolar_Hardy_Finsler} reduces to the following expression:
\begin{align} \label{eq1}
\int_M  |\nabla_g u|^2 dv_g  \geq 
& - \frac{(n-2)^2}{m^2} \int_M \bigg|  \sum_{i=1}^m \frac{\nabla_g d_i}{d_i} \bigg|^2  u^2 dv_g  \nonumber \\
& + \frac{n-2}{m} \sum_{i=1}^m  \int_M \mathrm{div} \Big(\frac{\nabla_g d_i}{d_i} \Big) u^2  dv_g.
\end{align}

Now we expand the first term of the right hand side. 
First of all, by using the eikonal equation \eqref{eq_eikonal}, one has
\begin{equation*}
\bigg| \frac{\nabla_g d_i}{d_i} - \frac{\nabla_g d_j}{d_j} \bigg|^2 = \frac{1}{d_i^2} + \frac{1}{d_j^2} - 2  \frac{g(\nabla_g d_i, \nabla_g d_j)}{d_i d_j} , 
\end{equation*}
for all $i,j \in \{1, \dots, m\}$. 

Then, using the 'expansion of the square' method and the eikonal equation again, we obtain 
\begin{align*} 
\bigg| \sum_{i=1}^m \frac{\nabla_g d_i}{d_i} \bigg|^2 & = \sum_{i,j=1}^m  g \left( \frac{\nabla_g d_i}{d_i}, \frac{\nabla_g d_j}{d_j} \right) \\
& =  \sum_{i=1}^m \frac{1}{d_i^2} + 2  \sum_{1 \leq i<j \leq m} \frac{g(\nabla_g d_i, \nabla_g d_j)}{d_i d_j}  \\
& =  m \sum_{i=1}^m \frac{1}{d_i^2} - \sum_{1 \leq i<j \leq m} \bigg| \frac{\nabla_g d_i}{d_i} - \frac{\nabla_g d_j}{d_j} \bigg|^2 .
\end{align*}

On the other hand, considering the second term of the right hand side of \eqref{eq1}, we have 
\begin{equation*}  
\mathrm{div} \Big(\frac{\nabla_g d_i}{d_i} \Big) = \frac{d_i \Delta_g d_i - 1}{d_i^2} , ~ \text{ for all } i=\overline{1,m}.
\end{equation*}

Substituting the expressions above and using direct calculations yields that 
\eqref{eq1}  is equivalent to the Riemannian multipolar inequality \eqref{multipolar_Hardy_curved}.
\end{remark}

Using Theorem \ref{thm_multipolar_Hardy_Finsler}, we can prove a bipolar Hardy inequality on complete Finsler manifolds.

\subsection*{\bf Proof of Theorem \ref{thm_bipolar_Hardy_Finsler}.}

Let $x_1, x_2 \in M$ be two distinct poles and $d_1, d_2: M \to [0, \infty)$ be the associated distance functions.
Using \eqref{proof1} and the eikonal equation \eqref{eq_eikonal}, we obtain
\begin{align*}  
&F^{*2} \Big(x, \frac{Dd_1}{d_1}(x) + \frac{Dd_2}{d_2}(x) \Big) \leq 
2 \Big( \frac{1}{d_1^2(x)} + \frac{1}{d_2^2(x)} \Big)   \\
&- l_F F^{*2}\Big(x, \frac{Dd_2}{d_2}(x) - \frac{Dd_1}{d_1}(x)\Big) , \text{ for a.e. } x \in M.
\end{align*}

Applying Theorem \ref{thm_multipolar_Hardy_Finsler} in the case $m=2$, then using the inequality above completes the proof of \eqref{bipolar_Hardy_Finsler}. \qed

\begin{remark}\rm 
Let $(B^n,F)$ be the usual Euclidean unit ball $B^n \subset \mathbb R^n$ endowed with the \textit{Funk metric} $F$, 
see Krist\'aly and Rudas \cite{KRu}. 
It turns out that $l_F=0$ and $r_F=+\infty$, thus both inequalities (\ref{multipolar_Hardy_Finsler})  and (\ref{bipolar_Hardy_Finsler})  reduce to  trivial statements. 
This particular example shows the importance of $l_F>0$ in Theorems \ref{thm_multipolar_Hardy_Finsler} and \ref{thm_bipolar_Hardy_Finsler}, respectively. 
\end{remark}

\section*{Acknowledgment} 
The authors are supported by the National Research, Deve\-lopment and Innovation Fund of Hungary, financed under the K$\_$18 funding scheme, Project No.  127926.

\end{document}